\documentclass{article}

\usepackage{graphicx}
\usepackage{psfrag}
\usepackage{amsmath}
\usepackage{amssymb}
\usepackage{epstopdf}
\usepackage{enumerate}
\usepackage{epstopdf}
\usepackage{amsthm}

\newtheorem{theorem}{Theorem}

\title{Torus knots are Fourier-$(1,1,2)$ knots}
\author{
  Jim Hoste\\
  Pitzer College}
\date{\today}
\ifx\pdfoutput\undefined
\else
\usepackage{graphicx,epstopdf}
\fi
\begin{document}
\maketitle

\begin{abstract}

Every torus knot can be represented as a Fourier-$(1,1,2)$ knot which is  the simplest possible Fourier representation for such a knot. This answers a question of Kauffman and confirms the  conjecture made by Boocher, Daigle, Hoste and Zheng.

In particular, the torus knot $T_{p,q}$ can be parameterized as
\begin{eqnarray*}
x(t)&=&\cos(pt)\\
y(t)&=&\cos(q t+\pi/(2 p))\\
z(t)&=&\cos(pt+\pi/2)+\cos((q-p)t+\pi/(2 p)-\pi/(4q)).
\end{eqnarray*}

\end{abstract}

\section{Introduction}

A Fourier-$(i,j,k)$ knot is one that can be parameterized as 
$$
\begin{array}{rcl}
     x(t) & = & A_{x,1}\cos(n_{x,1}t+\phi_{x,1})+...+A_{x,i}\cos(n_{x,i}t+\phi_{x,i}) \\
     y(t) & = & A_{y,1}\cos(n_{y,1}t+\phi_{y,1})+...+A_{y,j}\cos(n_{y,j}t+\phi_{y,j}) \\
     z(t) & = & A_{z,1}\cos(n_{z,1}t+\phi_{z,1})+...+A_{z,k}\cos(n_{z,k}t+\phi_{z,k}). \\
\end{array}
$$
Here the {\it frequencies}  $n_{x,r}, n_{y,r}$ and $n_{z,r}$ are integers and the {\it phase shifts} $\phi_{x,r}, \phi_{y,r}$ and $\phi_{z,r}$ and {\it amplitudes}  $A_{x,r}, A_{y,r}$ and $A_{z,r}$ are real numbers. 
Because any function can be closely approximated by a sum of cosines, every knot is a Fourier knot for some $(i,j,k)$. But a remarkable theorem of Lamm \cite{Lamm} states that in fact every knot is a Fourier-$(1,1,k)$ knot for some $k$.  The special case of Fourier-$(1,1,1)$ knots are called {\it Lissajous} knots. Lissajous knots are very highly symmetric: either strongly $+$~amphicheiral (if all three frequencies are odd), or 2-periodic and linking the axis of rotation once (if one frequency is even).  Hence most knots, including all torus knots,  are not Lissajous knots. Most of what is known about the topic can be found in   \cite{BHJ1994}--\cite{Lamm1996}.

While every knot can be expressed as a Fourier-$(1,1,k)$ knot for some integer $k$, it is not clear how big $k$ needs to be. To date no one has a produced a knot that requires $k$ to be bigger than 2. Could every knot be a Fourier-$(1,1,k)$  knot with $k\le 2$? It is known that this is true for all twist knots and for all 2-bridge knots to 14 crossings \cite{BDHZ2007}, \cite{HZ2006}. In this paper it is shown to be true for all torus knots. This was conjectured to be true in \cite{BDHZ2007} and answers the question posed by Kauffman in \cite{Kauffman}.

The parameterization of the torus knot $T_{p,q}$ given in Theorem~\ref{main theorem}  was found after the methodical sampling of Fourier-$(1,1,2)$ knots undertaken in \cite{BDHZ2007} yielded all torus knots to 16 crossings. Careful examination of these examples led to the general formula. 

\section{Torus Knots}

The standard parameterization of the $(p,q)$ torus knot as it winds round  an unknotted torus in $\mathbb R^3$, $p$ times in the longitudinal direction and $q$ times in the meridional direction, is
\begin{eqnarray*}
x(t)&=&R\cos(pt)+r \cos(pt)\cos(qt)\\
y(t)&=&R\sin(pt)+r \sin(pt)\cos(qt)\\
z(t)&=&r\sin(qt)
\end{eqnarray*}
Here the torus is defined by sweeping the circle of radius $r$ centered at the point $(R,0,0)$ and lying in the $xz$-plane around the $z$-axis. By using the identity
\begin{equation}
\cos x-\cos y=-2 \sin \frac{x+y}{2} \sin \frac{x-y}{2}
\label{trig identity}
\end{equation}
and converting sines to cosines, we may rewrite the parameterization as 
\begin{eqnarray*}
x(t)&=&R\cos(pt)+\frac{r}{2} \cos((p+q)t)+\frac{r}{2}\cos((p-q)t)\\
y(t)&=&R\cos(pt-\pi/2)+\frac{r}{2} \cos((p+q)t-\pi/2)+\frac{r}{2}\cos((p-q)t-\pi/2)\\
z(t)&=&r\cos(qt-\pi/2)
\end{eqnarray*}
Thus every torus knot can be expressed as Fourier-$(1,3,3)$ knot. This parameterization is given by both Costa \cite{Costa1990} and Kauffman  \cite{Kauffman}, and  Kauffman  asks if a simpler Fourier representation can be found. Since torus knots cannot be Lissajous, the following theorem is best possible.
\begin{theorem}
The torus knot $T_{p.q}$, with $0<p<q$, is equivalent to the Fourier-$(1,1,2)$ knot given by
\begin{eqnarray}
x(t)&=&\cos(pt)\label{xparam}\\
y(t)&=&\cos\left(q t+\pi/({2 p})\right)\label{yparam}\\
z(t)&=&\cos\left(p t+\pi/2\right )+\cos\left ((q-p)t+\pi/(2p)-\pi/(4 q)\right)\label{zparam}.
\end{eqnarray}
Furthermore, if $p$ is even, we may replace $\phi_{z,2}$ with $\pi /(2p)$.
\label{main theorem}
\end{theorem}

\begin{figure}[h]
    \begin{center}
    \leavevmode
    \scalebox{.9}{\includegraphics{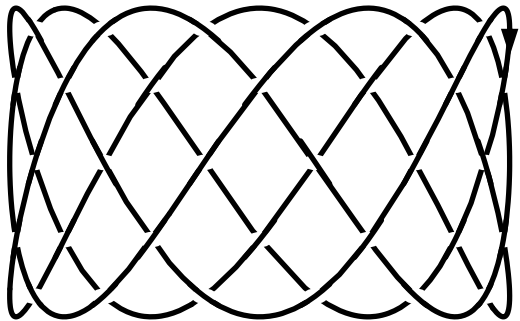}}
    \end{center}
\caption{The torus knot $T_{3,7}$ realized as a Fourier-$(1,1,2)$ knot.}
\label{37}
\end{figure}

\noindent{\bf Proof:} The $xy$-projection of the Fourier-$(1,1,2)$ knot given by Equations~\ref{xparam}-\ref{zparam} with $p=3$ and $q=7$ is shown in Figure~\ref{37}. The knot is oriented in the direction of increasing $t$ with the arrowhead placed at $t=0$. To see that this is indeed a diagram of the torus knot $T_{3,7}$ consider only those arcs that move from left to right. Notice that each of these arcs lies completely above any arc that is moving from right to left. Thus the right-going arcs form a braid that lies entirely above the braid formed by the left-going arcs. The two braids join to form a closed braid that we may think of as braided around the y-axis. This figure is representative of the general case and it will suffices to prove the following two claims:
\begin{enumerate}
\item Any string moving to the right  always crosses over any string moving to the left.
\item Any crossing between strings which are both moving to the right, or both to the left, is left handed.
\end{enumerate}
Suppose $t_1$ and $t_2$ are a pair of times that produce a double point in the projection of the knot into the $xy$-plane.  It is not hard to determine that the pair $(t_1, t_2)$ must be one of the following two types. The first possibility is that
\begin{equation*}
(t_1,t_2)=\left(  \frac{-k\pi}{p}+\frac{j\pi}{q}-\frac{\pi}{2 p q},  \frac{k\pi}{p}+\frac{j\pi}{q}-\frac{\pi}{2 p q}\right)
\end{equation*}
where $0<k<p$ and $1+\lfloor \frac{kq}{p}+\frac{1}{2p}\rfloor \le j \le \lfloor 2q- \frac{kq}{p}+\frac{1}{2p}\rfloor$. We call these {\it Type I} crossings. 

Or we may have
\begin{equation*}
(t_1,t_2)=\left(  \frac{-k\pi}{q}+\frac{j\pi}{p},  \frac{k\pi}{q}+\frac{j\pi}{p}\right)
\end{equation*}
where $0<k<q$ and $1+\lfloor \frac{pk}{q}   \rfloor \leq j \leq \lfloor 2p - \frac{pk}{q}   \rfloor$. We call these {\it Type II} crossings.

There are $pq-q$ crossings  of Type I, and $pq-p$ crossings of Type II.
See \cite{BDHZ2007} or \cite{HZ2006} for more details.

 At either type crossing, both arcs are moving to the right, or both to the left, if and only if $x'(t_1)x'(t_2)>0$. Evaluating $x'(t_1)x'(t_2)$ at a Type I crossing we obtain
 \begin{eqnarray*}
 x'(t_1)x'(t_2)&=&\left [ -p \sin\left(  -k\pi+\frac{j p\pi}{q}-\frac{\pi}{2  q}\right)\right]
 \left [ -p \sin\left(  k\pi+\frac{j p\pi}{q}-\frac{\pi}{2  q}\right)\right]\\
 &=&p^2\sin^2\left(  \frac{j p\pi}{q}-\frac{\pi}{2  q}\right)>0.
 \end{eqnarray*}
 Thus at a Type I crossing both strings are going to the right, or both are going to the left.
 A similar calculation reveals that at Type II crossings the strings are headed in opposite left-right directions. 
 
We may compute the sign of a crossing $(t_1, t_2)$ as follows. If $v(t)$ is the projection of the velocity vector into the $xy$-plane, then the cross product $v(t_1)\times v(t_2)$ will point up or down. Its dot product with the vector $(0,0,z(t_1)-z(t_2))$ will be positive if and only if the crossing is right handed. Thus the crossing has the same sign as 
\begin{equation}
\left[x'(t_1)y'(t_2)-x'(t_2)y'(t_1)\right]\left[z(t_1)-z(t_2)\right].
\label{crossing sign formula}
\end{equation}
Using Equation~\ref{trig identity} we may simplify $z(t_1)-z(t_2)$ as follows:
\begin{eqnarray}
z(t_1)-z(t_2)&=&\cos(pt_1+\pi/2)+\cos((q-p)t_1+\pi/(2p)-\pi/(4q))\nonumber\\
&&-\cos(pt_2+\pi/2)-\cos((q-p)t_2+\pi/(2p)-\pi/(4q))\nonumber \\
&=&-2 \sin \left (p\frac{t_1+t_2}{2}+\frac{\pi}{2}\right )\sin\left (p\frac{t_1-t_2}{2}\right)\label{z1-z2}\\
&&-2 \sin\left ((q-p)\frac{t_1+t_2}{2}+\frac{\pi}{2p}-\frac{\pi}{4q}\right )\sin\left ((q-p)\frac{t_1-t_2}{2}\right)\nonumber
\end{eqnarray}

To prove claim 1 we want to show that at every Type II crossing it is the upper string that is moving to the right. This is equivalent to showing that $x'(t_1)\left(z(t_1)-z(t_2)\right)>0.$ Using \ref{z1-z2} we obtain
\begin{eqnarray*}
x'(t_1)\left(z(t_1)-z(t_2)\right)&=&-p\sin(j\pi-pk\pi/q) [-2\sin(j \pi+\pi/2)\sin(-pk\pi/q) \\
&&\quad -2\sin((q-p)j\pi/p+\pi/(2p)-\pi/(4q))\sin((q-p)(-k\pi)/q) ]\\
&=& 2 p \sin^2(pk\pi/q)\left[ 1-(-1)^k\sin(jq\pi/p+\pi/(2p)-\pi/(4q))\right].
\end{eqnarray*}
This quantity is always positive since within the final factor the sine term has magnitude less than 1. Thus even if we are subtracting rather than adding, the difference is positive.

To prove claim 2, we must show that every Type I crossing is left handed. Computing the expression in  \ref{crossing sign formula}, one finds that its sign is the same as the sign of
$$-\sin(pj\pi/q-\pi/(2q))\sin(pj\pi/q-\pi/(4q))$$
which we want to show is negative. Notice that the two angles, $pj\pi/q-\pi/(2q)$ and $pj\pi/q-\pi/(4q)$ are very close together. In fact, they can never be separated by a multiple of $\pi$. For suppose that 
$$pj\pi/q-\pi/(2q)<r \pi<pj\pi/q-\pi/(4q)$$ 
for some integer $r$. Then, multiplying by $q/\pi$ gives
$$pj-1/2<r<pj-1/4$$
which is clearly impossible. Thus the sines of the two angles always have the same sign and hence their product is positive.

Claims 1 and 2 coupled with the fact that each string oscillates between 1 and $-1$ in the $y$-direction as it moves to the left or the right completes the proof that the closed braid is indeed a torus knot. Since there are a total of $pq-q$ Type I crossings, we see that it must be the knot $T_{p,q}$. 

The final assertion, that if $p$ is even the phase shift $\pi/(2p)-\pi/(4q)$ may be replaced with the simpler expression $\pi/(2p)$, can be proven by simply repeating the proof with the new phase shift. An alternative proof, and one that offers more insight into the situation, is to consider  the {\it phase torus}\\
 $\left \{ (\phi_{z,1}, \phi_{z,2})\mbox{ } | \mbox{ } 0\le \phi_{z,1}\le 2 \pi, 0\le \phi_{z,2}\le 2 \pi \right \}$. The Type I crossing with indices $(k,j)$ becomes singular on the horizontal line
 $$\phi_{z,2}=jp\pi/q+(1/p-1/q)/2\pi+m \pi$$
 where $m$ is an arbitrary integer. The Type II crossing with indices $(k,j)$ becomes singular along the diagonal lines
 $$\phi_{z,2}=\left \{ 
 \begin{array}{ll}
 (-1)^m \phi_{z,1}-jq\pi/p+m \pi, \mbox{ if $k$ is even}\\
 (-1)^{m+1}\phi_{z,1}-jq\pi/p+m \pi,  \mbox{ if $k$ is odd}
 \end{array}
 \right.
 $$
 where $m$ is an arbitrary integer.
 
 More detailed information on the phase torus and the kinds of singular curves that are associated to the Type I and II crossings can be found in \cite{BDHZ2007}. 
  Each region complementary to the singular lines defines a knot type. If $p$ is even, the two points $(\pi/2, \pi/(2p)-\pi/(4q))$ and $(\pi/2, \pi/(2p))$ lie in the interior of the same region and therefore define the same knot type. If $p$ is odd, the second point no longer lies in the same complementary region as the first point. In fact, it lies at the intersection of two  singular lines associated to Type II crossings and on the boundary of the region containing the first point.
\hfill$\square$

\pagebreak


\begin{thebibliography}{1}

\bibitem{BHJ1994}
M.~G.~V. Bogle, J.~E. Hearst, V.~F.~R. Jones, and L.~Stoilov.
\newblock Lissajous knots.
\newblock {\em J. Knot Theory Ramifications}, 3(2):121--140, 1994.

\bibitem{BDHZ2007}
Adam Boocher, Jay Daigle, Jim Hoste, and Wenjing Zheng.
\newblock Sampling Lissajous and Fourier knots.
\newblock {\em arXiv:0707.4210}, 2007.

\bibitem{HZ2006}
Jim Hoste and Laura Zirbel.
\newblock Lissajous knots and knots with Lissajous projections.
\newblock {\em arXiv: math.GT/0605632}, 2006.

\bibitem{JP1998}
Vaughan F.~R. Jones and J\'{o}zef~H. Przytycki.
\newblock Lissajous knots and billiard knots.
\newblock {\em Banach Center Publications}, (42):145--163, 1998.

\bibitem{Kauffman}
Louis Kauffman.
\newblock Fourier knots.
\newblock {\em arXiv: q-alg/9711013}.

\bibitem{Lamm}
Christoph Lamm.
\newblock Fourier knots.
\newblock {\em Preprint}.

\bibitem{Lamm1996}
Christoph Lamm.
\newblock There are infinitely many lissajous knots.
\newblock {\em Manuscripta Math.}, 93:29--37, 1996.

\bibitem{Costa1990}
Sueli~I. Rodrigues~Costa.
\newblock On closed twisted curves.
\newblock {\em Proc. Amer. Math. Soc.}, 109(1):205--214, 1990.

\end{thebibliography}
\end{document}